\documentclass[12pt]{amsart}

\newtheorem{teo}{Theorem}[section]
\newtheorem{lema}{Lemma}[section]

\def\dem{\noindent{\bf Proof: }}
\def\ve{\varepsilon}
\def\RR{{\mathbb{R}}}

\newcommand{\pf}{\noindent{\bf{Proof:  }}}

\begin{document}

\title [The blow-up problem]{The blow-up
problem for a semilinear
parabolic equation with a potential}

\author[C. Cortazar\and M. Elgueta \and J.D. Rossi]{Carmen Cortazar
\and Manuel Elgueta \and Julio D. Rossi}
\thanks{Supported by Universidad de Buenos Aires under grant TX048,
by ANPCyT PICT No. 03-00000-00137 and CONICET (Argentina) and by
Fondecyt 1030798 and Fondecyt Coop. Int. 7050118 (Chile).
\newline
\noindent 2000 {\it Mathematics Subject Classification }
35K57, 35B40.}
\keywords{Blow-up, semilinear parabolic equations.}
\address{Departamento  de Matem{\'a}tica, Universidad Catolica de Chile,
\hfill\break\indent Casilla 306, Correo 22, Santiago, Chile. }
\email{\tt ccortaza@mat.puc.cl, melgueta@mat.puc.cl}

\address{
Instituto de Matem{\'a}ticas y F\'{\i}sica Fundamental
\hfill\break\indent Consejo Superior de Investigaciones
Cient\'{\i}ficas \hfill\break\indent Serrano 123, Madrid,
Spain,\hfill\break\indent on leave from Departamento de
Matem{\'a}tica, FCEyN UBA (1428)\hfill\break\indent Buenos Aires,
Argentina.} \email{{\tt jrossi@dm.uba.ar }}


\begin{abstract} Let $\Omega$ be a bounded smooth domain in $\RR^N$.
We consider the problem $u_t= \Delta u + V(x) u^p$ in $\Omega
\times [0,T)$, with Dirichlet boundary conditions $u=0$ on
$\partial \Omega \times [0,T)$ and initial datum $u(x,0)= M
\varphi (x)$ where $M
\geq 0$, $\varphi$ is positive and compatible with the boundary
condition. We give estimates for the blow up time of solutions for
large values of $M$. As a consequence of these estimates we find
that, for $M$ large, the blow up set concentrates near the points
where $\varphi^{p-1}V$ attains its maximum.
\end{abstract}

\maketitle

\bigskip

\section*{{\bf 1. INTRODUCTION}}
\setcounter{section}{1} \setcounter{equation}{0}

\bigskip

In this paper we study the blow-up phenomena for the following
semilinear parabolic problem with a potential
\begin{equation}
\label{problema}
\begin{array}{ll}
u_t = \Delta u + V(x) \, u^p \qquad & \mbox{ in }\Omega \times
(0,T),\\[8pt]
u(x,t) = 0  \qquad & \mbox{ on } \partial \Omega \times (0,T), \\[8pt]
u(x,0) = M\varphi(x)  \qquad & \mbox{ in } \Omega .
\end{array}
\end{equation}

First, let us state our basic assumptions. They are: $\Omega$ is a
bounded, convex, smooth domain in $\RR^N$ and the exponent $p$ is
subcritical, that is, $1 < p < (N+2)/(N-2)$. The potential $V$ is
Lipschitz continuous and there exists a constant $c>0$ such that
$V(x) \geq c$ for all $x
\in
\Omega$. As for the initial condition we assume that $M \geq 0$
and that $\varphi$ is a smooth positive function compatible with
the boundary condition. Moreover, we impose that
\begin{eqnarray}\label{initial.2} M\Delta
\varphi + \frac{\min\limits_{x \in \Omega} V(x)}{2} M^p \varphi^p
\ge 0.
\end{eqnarray}
We note that \eqref{initial.2} holds for $M$ large if $\Delta
\varphi$ is nonnegative in a neighborhood of the set where $\varphi$
vanishes.

It is known that, and we will prove it later for the sake of
completeness, once $\varphi$ is fixed the solution to
(\ref{problema}) blows up in finite time for any $M$ sufficiently
large. By this we understand that there exists a time $T=T(M)$
such that $u$ is defined in $\Omega \times [0,T)$ and
$$\lim_{t\to T} \|
u(\cdot,t)\|_{L^\infty(\Omega)}=+\infty . $$

The study of the blow-up phenomena for parabolic equations and
systems has attracted considerable attention in recent years, see
for example, \cite{B}, \cite{BB}, \cite{GK1}, \cite{GK2},
\cite{GV}, \cite{HV1}, \cite{HV2}, \cite{M}, \cite{Z} and the
corresponding references. A good review in the topic can be found
in \cite{GaVa}. When a large or small diffusion is considered, see
\cite{IY}, \cite{MY}.

Important issues in a blow-up problem are to obtain estimates for
the {\it blow-up time}, $T(M)$, and determine the spatial
structure of the set where the solution becomes unbounded, that
is, the {\it blow-up set}. More precisely, the blow-up set of a
solution $u$ that blows up at time $T$ is defined as
$$
B(u) = \{ x / \ \mbox{there exist }x_n \to x, \, t_n\nearrow T, \
\mbox{with } u(x_n, t_n) \to \infty \}.
$$
The problem of estimating the blow-up time and the description and
location of the blow-up set has proved to be a subtle problem and
has been addressed by several authors. See for example
\cite{SGKM}, \cite{GaVa} and the corresponding bibliographies.

\medskip

Our interest here is the description of the asymptotic behavior of
the blow-up time, $T(M)$, and of the blow-up set, $B(u)$, as $M
\to \infty$. It turns out that their asymptotics depend on a combination of
the shape of both the initial condition, $\varphi$, and the
potential $V$. Roughly speaking one expects that if $\varphi
\equiv 1$ then the blow-up set should concentrate near the points
where $V$ attains its maximum. On the other hand if $V\equiv 1$
the blow-up set should be near the points where $\varphi$ attains
its maximum. Just to see what to expect, if we drop the laplacian,
we get the ODE $ u_t = V(x) u^p$ with initial condition $u(x,0) =
M \varphi (x)$. Here $x$ plays the role of a parameter. Direct
integration gives $ u(x,t) = C (T-t)^{-1/(p-1)}$ with
$$T=\frac{ M^{1-p}}{ (p-1) V(x)
\varphi^{p-1}(x) }.$$ Hence, blow-up takes place at points $x_0$ that
satisfy $ V(x_0) \varphi^{p-1} (x_0) = \max_x V (x)
\varphi^{p-1}$. Therefore, we expect that the quantity that plays
a major role is $(\max_x V (x) \varphi^{p-1}(x) )$.

\begin{teo} \label{teo.Mej}
There exists $\bar M >0$ such that if $M \geq \bar M$ the solution
of \eqref{problema} blows up in a finite time that we denote by
$T(M)$. Moreover, let
$$
A= A(\varphi, V) := \frac{1}{ (\max_x \varphi^{p-1}(x) V (x))  },
$$
then there exist two positive constants $C_1$, $C_2$, such that,
for $M$ large enough,
\begin{equation}
\label{T.2}
\displaystyle  - \frac{C_1}{
M^{\frac{p-1}{4} }  }
\le T(M) M^{p-1} - \frac{A}{p-1} \displaystyle  \le
\frac{C_2}{ M^{\frac{p-1}{3} }  },
\end{equation}and the blow-up set verifies,
\begin{equation}\label{V(BU).2}
\varphi^{p-1} (a) V (a) \ge \frac{1}{A} -\frac{C}{ M^{\gamma}  },
\qquad
\mbox{ for all } a \in B(u),
\end{equation}
where  $\gamma= min (\frac{p-1}{4},\frac{1}{3})$.
\end{teo}

Note that this result implies that
$$
\lim_{M\to \infty} T(M) M^{p-1} =  \frac{A}{p-1}.
$$
Moreover, it provides precise lower and upper bounds on the
difference $T(M) M^{p-1} -  \frac{A}{p-1}$.

We also observe that \eqref{V(BU).2} shows that the set of blow-up
points concentrates for large $M$ near the set where
$\varphi^{p-1} V$ attains its maximum.

If in addition the potential $V$ and the initial datum $\varphi$
are such that $\varphi^{p-1} V$ has a unique non degenerate
maximum at a point $\bar a$, then there exist constants $c >0$ and
$d>0$ such that
$$ \varphi^{p-1} (\bar a) V(\bar a) -  \varphi^{p-1} (x) V(x)
\ge c |\bar a - x |^2 \quad \mbox{ for all } x \in B(\bar a,d).$$
Therefore, according to our result, if $M$ is large enough one has
$$
|\bar a -a| \le \frac{C}{ M^{\frac{\gamma}{2}}} \quad \mbox{ for
any } a\in B(u),
$$
with $\gamma= min (\frac{p-1}{4},\frac{1}{3})$.

\medskip

Throughout the paper we will denote by $C$ a constant that does
not depends on the relevant parameters involved but may change at
each step.


\section*{\bf {2. Proof of Theorem \ref{teo.Mej}.}}
\setcounter{equation}{0} \setcounter{section}{2} We begin with a
lemma that provides us with an upper estimate of the blow-up time.
This upper estimate gives the upper bound for $T(M) M^{p-1}$ in
\eqref{T.2} and will be crucial in the rest of the proof of Theorem
\ref{teo.Mej}.

\begin{lema}\label{est.T} There exist a constant $C>0$ and $M_0>0$ such
that for every $M\ge M_0$, the solution of \eqref{problema} blows
up in a finite time that verifies
\begin{equation}\label{esti.T}
T (M) \le  \frac{A}{M^{p-1}
(p-1)}+\frac{C}{M^{\frac{p-1}{3}}M^{p-1}}.
\end{equation}
\end{lema}

\pf Let $\bar a \in \Omega$ be such that
$$\varphi^{p-1} (\bar a)V(\bar
a)=\max_x \varphi^{p-1}(x)V(x),
$$
$L$ the constant of Lipschitz continuity of $V$, and $K$ an upper
bound for the first derivatives of $\varphi$ and $L$.

In order to get the upper estimate let $M$ be fixed and
$\ve=\ve(M)
>0$ to be defined latter, small enough so all functions involved are well defined.
Pick
$$
\delta =\frac{\ve}{2 K},
$$
then
$$
 V(x) \ge  V(\bar a)
 -\frac{\ve}{2}  \quad
 \mbox{ and }  \qquad \varphi(x) \ge \varphi(\bar a)- \ve
 \quad \mbox{ for all } x \in B(\bar a,
 \delta).
$$
Let $w$ be the solution of
$$\begin{array}{ll}
 w_t = \Delta w + \left(V(\bar a) -\displaystyle
\frac{\ve}{2}\right) w^p & \mbox{ in }
 B(\bar a,\delta) \times (0,T_w), \\[8pt]
w=0 & \mbox{ on }  \partial B(\bar a, \delta) \times (0,T_w),
\\[8pt]
w(x,0) = M (\varphi (\bar a)- \ve), & \mbox{ in }
  B(\bar a,\delta)
\end{array}$$
and $T_w$ its corresponding blow up time. A comparison argument
shows that $u \ge w$ in $B(\bar a, \delta) \times (0,T)$ and hence
$$ T \le T_w. $$

Our task now is to estimate $T_w$ for large values of $M$. To this
end, let $\lambda_1 (\delta)$ be the first eigenvalue of $-\Delta$
in $B(\bar a,\delta)$ and let $\phi_1$ be the corresponding
positive eigenfunction normalized so that
$$
\int_{B(\bar a, \delta) }
\phi_1 (x) \, dx =1.
$$
That is,
$$
\left\{\begin{array}{ll} - \Delta \phi_1 = \lambda_1
(\delta) \phi_1  , \qquad &
\mbox{ in } B(\bar
a,\delta),\\[8pt]
\phi_1 =0 \qquad &  \mbox{ on } \partial B(\bar a, \delta).
\end{array}\right.
$$
Now, set
$$
\Phi
(t) = \int_{B(\bar a,\delta)} w(x,t) \phi_1 (x)\, dx.
$$
Then $\Phi(t)$ satisfies $\Phi(0) = M(\varphi (\bar a)- \ve)$ and
$$
\begin{array}{rl}
\Phi ' (t) & = \displaystyle \int_{B(\bar a,\delta)} w_t (x,t) \phi_1 (x) \, dx \\[10pt]
& = \displaystyle  \int_{B(\bar a,\delta)} \left( \Delta w(x,t)
\phi_1 (x)  + \left( V(x_1) -\frac{\ve}{2} \right)  w^p (x,t)
\phi_1(x) \right) \, dx \\[10pt]
& \ge \displaystyle -\lambda_1 (\delta) \int_{B(\bar a,\delta)}
w(x,t)
\phi_1 (x) \, dx \\[10pt]
& \qquad \displaystyle+ \left( V(\bar a) -\frac{\ve}{2} \right)
\left(
\int_{B(\bar a,\delta)} w (x,t) \phi_1(x) \, dx
 \right)^p \\[10pt]
& = \displaystyle -\lambda_1 (\delta) \Phi (t) + \left( V(\bar a)
-
\frac{\ve}{2}\right ) \Phi (t)^p.
\end{array}
$$

Let us recall that there exists a constant $D$, depending on the
dimension only, such that the eigenvalues of the laplacian scale
according to the rule $
\lambda_1 (\delta ) = D
\delta^{-2}$.

Now, we choose $\ve$ such that
$$\lambda_1 (\delta) = D \delta^{-2}=D \left(\frac{\ve}{2 K}\right)^{-2}=
\frac{\ve}{2}(M(\varphi (\bar a)- \ve))^{p-1}.$$
So, $\ve$ is of order
$$
\ve \sim \frac {C}{ M^{\frac{p-1}{3}}}.
$$
Choose $M_0$ such that for $M\geq  M_0$  the resulting $\ve$ is
small enough. Then for any $M \geq M_0$ we have that
\begin{equation}\label{ec.phi}
\Phi ' (t)
\geq (V(\bar a) - \ve) \Phi (t)^p,
\end{equation}
for all $t \geq 0$ for which $\Phi$ is defined.

Since $\Phi (0)=M (\varphi(\bar a) - \ve)$ and $T_w$ is less or
equal than the blow up time of $\Phi$ integrating \eqref{ec.phi}
it follows that
$$
\begin{array}{rl}
T_w & \le \displaystyle
\frac{1}{M^{p-1} (p-1)  ( V(\bar a) -\ve ) (\varphi (\bar a)
-\ve)^{p-1} } \\[12pt]
& \displaystyle \leq \frac{1}{M^{p-1} (p-1)   V(\bar a) \varphi
(\bar a) ^{p-1}}+\frac{C}{M^{\frac{p-1}{3}}M^{p-1} },
\end{array}
$$
for all $M \geq M_0$. \qed

\medskip

Now we prove a lemma that provides us with an upper bound for the
blow up rate. We observe that this is the only place where we use
hypothesis \eqref{initial.2}.

\begin{lema}\label{cota.por.arriba}
Assume \eqref{initial.2}. Then there exists a constant $C$
independent of $M$ such that
$$
u(x,t) \le C(T-t)^{-\frac{1}{p-1}}.
$$
\end{lema}

\dem Let $m = \min\limits_{x \in \Omega} V$. Following ideas of
\cite{FMc}, set
$$ v = u_t-\frac{m}{2} u^p.$$ Then $v$ verifies
$$ \begin{array}{ll}
v_t-  \Delta v  - V(x)p u^{p-1}v =
\displaystyle \frac{m}{2} p(p-1) u^{p-2}
|\nabla u |^2 \geq 0  & \mbox{ in }\Omega \times
(0,T),\\[8pt]
v =0   & \mbox{ on } \partial \Omega \times (0,T), \\[8pt]
v(x,0) = M\Delta \varphi + \left(V(x)-\displaystyle
\frac{m}{2}\right) M^p \varphi^p
\geq 0
 &\mbox{ in } \Omega.
\end{array}
$$
Therefore $v \geq 0$ and hence
$$u_t \geq\frac{m}{2} u^p.$$
Integrating this inequality from $0$ to $T$ we get
$$ u(x,t)\leq \frac {2^{\frac{1}{p-1}}}{( m (p-1)
(T-t)))^{\frac{1}{p-1}}}\equiv C(T-t)^{-\frac{1}{p-1}},$$ as we
wanted to prove.
\qed

\medskip

We are now in a position to prove Theorem \ref{teo.Mej}.

\medskip

\noindent{\bf Proof of Theorem \ref{teo.Mej}:} The idea of the
proof is to combine the estimate of the blow-up time proved in
Lemma \ref{est.T} with local energy estimates near a blow-up point
$a$, like the ones considered in \cite{GK1} and \cite{GK2}, to
obtain an inequality that forces $\varphi^{p-1} (a)V(a)$ to be
close to $\max_x \varphi^{p-1}V$.

\medskip

Let us now proceed with the proof of the estimates on the blow-up
set. We fix for the moment $M$ large enough such that $u$ blows up
in finite time $T=T(M)$ and let $a=a(M)$ be a blow up point. As in
\cite{GK2}, for this fixed $a$ we define
$$
w(y, s) = (T-t)^{\frac{1}{p-1}} u(a+y(T-t)^{\frac{1}{2}},t) |_{t=
T(1-e^{-s})} .
$$
Then $w$ satisfies
\begin{equation} \label{ec.w}
w_s =  \Delta w - \frac{1}{2} y  \cdot \nabla w -\frac{1} { p-1}w
+ V(a+yTe^{-\frac{s}{2}}) w^{p},
\end{equation}
in $\cup_{s \in (0,\infty)} \Omega (s) \times  \{s\}$ where
$\Omega (s) =
\Omega_a(s) =
\{ y
\, : \, a+yTe^{-\frac{s}{2}} \in \Omega \}$ with $w(y,0) =
T^{\frac{1}{p-1}} \varphi(a+y T^{\frac{1}{2}})$. The above
equation can rewritten as
$$
w_s =  \frac{1} { \rho} \nabla (\rho \nabla w)-  \frac{1} { p-1}w
+ V(a+y T e^{-\frac{s}{2}}) w^{p}
$$
where $\rho(y)=\exp(\frac{-|y|^2} { 4})$.

Consider the energy associated with the ''frozen'' potential
$$
V
\equiv V(a),
$$
that is
$$
E(w)=\int_{\Omega(s)} \left({\frac{1}{2}} |\nabla w|^2
+{\frac{1}{2(p-1)}} w^2 -{\frac{1}{p+1}} V(a) w^{p +1}
\right)\rho(y)\, dy .$$ Then, using the fact that $\Omega$ is
convex, we get
$$\frac{dE}{ds}\le -\int_{\Omega(s)}(w_s)^2 \rho(y)\, dy + \int_{\Omega(s)}
(V(a+yTe^{-\frac{s}{2}})-V(a)) w^p w_s \rho(y)\, dy. $$

Since $V(x)$ is Lipschitz and $w$ is bounded due to Lemma
\ref{cota.por.arriba}, then there exists a constant $C$ depending
only on $N$, $p$ and $V$, recall that the constant in Lemma
\ref{cota.por.arriba} does not depend on $M$, such that
$$\frac{dE}{ds} \leq -\int (w_s)^2 \rho(y)\, dy +C
e^{-\frac{s}{2}}T \left(\int (w_s)^2 \rho(y)\, dy\right)^{1/2}.$$

Maximizing the right hand side of the above expression with
respect to $\int (w_s)^2 \rho(y)\, dy$ we obtain
$$\frac{dE}{ds} \leq C
e^{-s }T^2
$$
and integrating is $s$ we get
\begin{equation}
\label{E.decrece}
E(w)\leq E(w_0)+ C T^2 .
\end{equation}

Since $w$ is bounded and satisfies \eqref{ec.w}, following the
arguments given in
\cite{GK1} and \cite{GK2}, one can prove that $w$ converges as
$s\to \infty$ to a non trivial bounded stationary solution of the
limit equation
\begin{equation}\label{ec.limite}
0=  \Delta z - \frac{1} { 2} y \cdot \nabla z -\frac{1} { p-1}z +
V(a) z^{p}
\end{equation}
in the whole $\RR^N$.

Again by the results of \cite{GK1} and
\cite{GK2}, since $p$ is subcritical, $1<p <(N+2)/(N-2)$,
the only non trivial bounded positive solution of
(\ref{ec.limite}) with $V(a) =1$ is the constant
$(p-1)^{-\frac{1}{p-1}}$. A scaling argument gives that the only
non trivial bounded positive solution of (\ref{ec.limite})
 is the constant $k
= k(a)$ given by
$$k(a) =\frac {1}{( V(a) (p-1))^{\frac{1}{p-1}}}.$$
Therefore, we conclude that
$$\lim_{s \rightarrow \infty}w=k(a)$$
if $a$ is a blow-up point. Also by the results of \cite{GK1},
\cite{GK2} we have
\begin{equation}\label{manuel}
E(w(\cdot ,s)) \rightarrow E(k(a)) \qquad \mbox{ as } s\to
\infty,
\end{equation}
where
$$
\begin{array}{rl}
E(k(a))= & \displaystyle
 \int \left({\frac{1}{2(p-1)}}(k(a))^2 - {\frac{1}{p+1}}V(a) (k(a))^{p +1}
 \right)\rho(y) \, dy \\
\\
= & (k(a))^2 \displaystyle
\left({\frac{1}{2(p-1)}}-{\frac{1}{(p+1)(p-1)}}\right) \int \rho(y) \, dy.
\end{array}
$$

By \eqref{E.decrece} and \eqref{manuel} we obtain that, if $a$ is
a blow-up point, then
$$
E(k(a))\leq E(w_0)+ C T^2 .
$$
where $w_0(y) = w(y,0) = T^{\frac{1}{p-1}}M \varphi(a + y
T^{\frac{1}{2}}).$

As $\varphi$ is smooth, $y \rho(y)$ integrable, and
$T^{\frac{1}{p-1}}M$ is bounded by Lemma \ref{est.T}, there are
constants $C$ independent of $a$ such that for $ M\geq M_0$
$$
\begin{array}{rl}
E(w(\cdot ,0)) &
\displaystyle =\int_{\Omega(0)} \left({\frac{1}{2}} |\nabla w_0 (y)|^2
+{\frac{1}{2(p-1)}} w_0^2 (y)
\right)\rho(y)\, dy \\[10pt]
& \qquad \displaystyle -\int_{\Omega(0)} \left({\frac{1}{p+1}}
V(a) w_0^{p +1} (y)
\right)\rho(y)\, dy \\[10pt]
& \displaystyle \leq\int_{\Omega(0)} \left({\frac{1}{2}}
(T^{\frac{1}{p-1}}M)^2 T|\nabla \varphi(a)|^2\right)\rho(y)\, dy
\\[10pt]
&\qquad \displaystyle +\int_{\Omega(0)}
\left({\frac{1}{2(p-1)}} (T^{\frac{1}{p-1}}M \varphi(a))^2  \right)\rho(y)\, dy \\[10pt]
& \qquad \displaystyle - \int_{\Omega(0)}
\left({\frac{1}{p+1}}
V(a) (T^{\frac{1}{p-1}}M \varphi(a))^{p +1} \right)\rho(y)\, dy \\[10pt]
& \qquad + CT^{\frac{3}{2}}+ CT^{\frac{1}{2}}.
\end{array}
$$
Therefore, since $|\nabla \varphi|$ is bounded,
$$
\begin{array}{rl}
E(w(\cdot ,0)) & \displaystyle \leq \int_{\Omega (0)}
\left({\frac{1}{2(p-1)}}
(T^{\frac{1}{p-1}}M \varphi(a))^2 \right)\rho(y)\, dy \\[10pt]
& \qquad - \displaystyle \int_{\Omega (0)} \left({\frac{1}{p+1}}
V(a) (T^{\frac{1}{p-1}}M \varphi(a))^{p +1}
\right)\rho(y)\, dy
\\[10pt]
& \qquad + CT^{\frac{3}{2}}+ CT^{\frac{1}{2}}.
\end{array}
$$
Or, since $T\le 1$ for $M$ large
$$
E(w(\cdot ,0))  \le E(T^{\frac{1}{p-1}}M \varphi(a))+
CT^{\frac{1}{2}}.
$$

Hence we arrive to the following bound for $E(k(a))$
\begin{equation} \label{carmen.*}
 E(k(a))\leq E(w(\cdot ,0))+ C T^2 \leq E(T^{\frac{1}{p-1}}M
\varphi(a))+ CT^{\frac{1}{2}}.
\end{equation}
Observe that if $b$ is a constant then the energy can be written
as
$$
E(b)=\Gamma F(b),
$$
where $\Gamma$ is the constant
$$\Gamma
=\int \rho (y) \, dy
$$
and $F$ is the function
$$F(z)=\left({\frac{1}{2(p-1)}} z^2 -{\frac{1}{p+1}} V(a) z^{p +1}
\right).$$

As $F$ attains a unique maximum at $k(a)$ and $F''(k(a))=-1$ there
are $\alpha$ and $\beta$ such that if $|z-k(a)|\leq \alpha $ then
$$F''(z)\leq -\frac{1}{2},$$ and if $|F(z)-F(k(a))|\leq \beta$ then
$$|z-k(a)|\leq \alpha.$$

{}From \eqref{carmen.*} we obtain
$$
F(k(a))\leq F(T^{\frac{1}{p-1}}M \varphi(a))+ CT^{\frac{1}{2}}.
$$
If $M_1$ is such that $C(T(M_1))^{\frac{1}{2}}= \beta$ then for
$M\geq  max(M_0,M_1)$
$$
\beta \ge CT^{\frac{1}{2}}\geq F(k(a))- F(T^{\frac{1}{p-1}}M \varphi(a)).
$$
Hence by the properties of $F$,
$$
|k(a)- T^{\frac{1}{p-1}}M \varphi(a)| \le \alpha.
$$
Therefore
$$
CT^{\frac{1}{2}}\geq F(k(a))- F(T^{\frac{1}{p-1}}M \varphi(a))\geq
\frac{1}{4} (T^{\frac{1}{p-1}}M \varphi(a)-k(a))^2.
$$
So, using Lemma \ref{est.T},
\begin{equation}\label{carmen.**}
\begin{array}{rl}
k(a)-CT^{\frac{1}{4}} & \displaystyle \leq T^{\frac{1}{p-1}}M
\varphi(a)
\\[8pt]
& \displaystyle \leq
\frac{\varphi(a)}{ (p-1)^{\frac{1}{p-1}}   V^{\frac{1}{p-1}}(\bar a) \varphi (\bar a)
}+\frac{C \varphi(a)}{M^{\frac{1}{3}}}
\\[12pt]
& \displaystyle= k(a)\theta(a)+\frac{C
\varphi(a)}{M^{\frac{1}{3}}},
\end{array}
\end{equation}
where
$$
\theta(a)=\left(\frac{\varphi(a)V(a)^{\frac{1}{p-1}}}{\varphi(\bar a )V(\bar
a)^{\frac{1}{p-1}}}\right)
$$
and $\bar a$ is such that
$$\varphi^{p-1} (\bar a)V(\bar
a)=\max_x \varphi^{p-1}(x)V(x).
$$
Recall that
$$T\leq \frac{C}{M^{p-1}}.$$
Therefore, we get
$$k(a)(1-\theta(a))\leq  \frac{C \varphi (a)}{M^{\frac{1}{3}}} + \frac{C}{M^{\frac{p-1}{4}}}
 \leq\frac{C}{M^{\gamma}},
 $$
with $\gamma= min (\frac{p-1}{4},\frac{1}{3})$.

As $V$ is bounded we have that $k(a)$ is bounded from below, hence
$$(1-\theta(a))\leq \frac{C}{M^{\gamma}},$$
that is,
$$
\theta(a)\geq 1-\frac{C}{M^{\gamma}}
$$
and we finally obtain
\begin{equation}\label{jorge}
\varphi(a)V(a)^{\frac{1}{p-1}}\geq \varphi(\bar a )V(\bar
a)^{\frac{1}{p-1}} -\frac{C}{ M^{\gamma}}.
\end{equation}
This proves \eqref{V(BU).2}.

\medskip

To obtain the lower estimate for the blow-up time observe that
from \eqref{jorge} and the fact that $V(a)\geq c>0$ we get
\begin{equation}\label{pepe}
\begin{array}{rl}
\varphi(a) & \displaystyle \geq\frac{\varphi(\bar a )V(\bar
a)^{\frac{1}{p-1}}}{V(a)^{\frac{1}{p-1}}}
-\frac{C}{V(a)^{\frac{1}{p-1}} M^{\gamma}} \\[10pt]
& \displaystyle
\geq\frac{\varphi(\bar a
)V(\bar a)^{\frac{1}{p-1}}}{V(a)^{\frac{1}{p-1}}} -\frac{C}{
M^{\gamma}}\\[10pt]
& \geq C>0.
\end{array}
\end{equation}

Inequality \eqref{carmen.**} gives us
 $$
\frac {1}{( V(a) (p-1))^{\frac{1}{p-1}}}- C T^{\frac{1}{4}}
 \leq T^{\frac{1}{p-1}}M
\varphi(a).
$$
Hence
$$
\frac {1}{ \varphi (a) ( V(a) (p-1))^{\frac{1}{p-1}}}-\frac{C T^{\frac{1}{4}}}{\varphi (a) }
 \leq T^{\frac{1}{p-1}}M.
$$
By \eqref{pepe} and $\varphi^{p-1} (\bar a)V(\bar a)=\max_x
\varphi^{p-1}(x)V(x)$ we get
$$
\frac {1}{ \varphi (\bar a) ( V(\bar a) (p-1))^{\frac{1}{p-1}}}- C T^{\frac{1}{4}}
 \leq T^{\frac{1}{p-1}}M
$$
and using
$$
T \le \frac{C}{M^{p-1}}
$$
we obtain
$$
\frac {1}{ \varphi (\bar a) ( V(\bar a) (p-1))^{\frac{1}{p-1}}}-
\frac{C}{M^{\frac{p-1}{4}}}
 \leq T^{\frac{1}{p-1}}M
$$
as we wanted to prove. \qed


\begin{thebibliography}{BHRP}


\bibitem[B]{B} J. Ball. {\it Remarks on blow-up and nonexistence theorems
for nonlinear evolution equations}. Quart. J. Math. Oxford, Vol.
28, (1977), 473--486.

\bibitem[BB]{BB} C. Bandle and H. Brunner. {\it Blow-up in diffusion
equations: a survey.} J. Comp. Appl. Math. Vol. 97, (1998), 3--22.

\bibitem[FMc]{FMc} A. Friedman and J. B. Mc Leod.
{\it Blow up of positive solutions of semilinear heat equations}.
Indiana Univ. Math. J., Vol. 34, (1985), 425--447.

\bibitem[GV]{GV} V. A. Galaktionov and J. L. Vazquez. {\it
Continuation of blow-up solutions of nonlinear heat equations in
several space dimensions}. Commun. Pure Applied Math. 50, (1997),
1--67.

\bibitem[GV2]{GaVa} V. A. Galaktionov and J. L. V{\'a}zquez.
{\it The problem of blow-up in nonlinear parabolic equations}.
Discrete Contin. Dynam. Systems A. Vol 8, (2002), 399--433.


\bibitem[GK1]{GK1} Y. Giga and R. V. Kohn. {\it Nondegeneracy of blow up
for semilinear heat equations}. Comm. Pure Appl. Math. Vol. 42,
(1989), 845--884.

\bibitem[GK2]{GK2} Y. Giga and R. V. Kohn. {\it Characterizing blow-up
 using similarity variables}. Indiana Univ. Math. J. Vol. 42,
(1987), 1--40.



\bibitem[HV1]{HV1} M. A. Herrero and J. J. L. Velazquez. {\it
Flat blow up in one-dimensional, semilinear parabolic problems},
Differential Integral Equations. Vol. 5(5), (1992), 973--997.

\bibitem[HV2]{HV2} M. A. Herrero and J. J. L. Velazquez. {\it
Generic behaviour of one-dimensional blow up patterns}. Ann.
Scuola Norm. Sup. di Pisa, Vol. XIX (3), (1992), 381--950.


\bibitem[IY]{IY} K. Ishige and H. Yagisita. {\it Blow-up
problems for a semilinear heat equation with large diffusion}. J.
Differential Equations. Vol. 212(1), (2005), 114--128.


\bibitem[M]{M} F. Merle. {\it Solution of a nonlinear heat equation with
arbitrarily given blow-up points}. Comm. Pure Appl. Math. Vol.
XLV, (1992), 263--300.

\bibitem[MY]{MY} N. Mizoguchi and E. Yanagida. {\it Life span of
solutions for a semilinear parabolic problem with small
diffusion}. J. Math. Anal. Appl. Vol. 261(1), (2001), 350--368.

\bibitem[SGKM]{SGKM} A. Samarski, V. A. Galaktionov, S. P. Kurdyunov
and A. P. Mikailov. Blow-up in quasilinear parabolic equations.
Walter de Gruyter, Berlin, (1995).

\bibitem[Z]{Z} H. Zaag. {\it One dimensional behavior of singular
$N$ dimensional solutions of semilinear heat equations}. Comm.
Math. Phys. Vol. 225 (3), (2002), 523--549.

\end{thebibliography}
\end{document}